\documentclass[a4paper,12pt, reqno]{amsart}
\usepackage{a4wide}
\usepackage{amsthm}
\usepackage{amssymb}
 \newtheorem{lemma}{Lemma}

 \theoremstyle{remark}

 \numberwithin{equation}{section}

\begin{document}

\title{A note on Hardy's theorem }
\author{Usha  K. Sangale}
\address{SRTM University,  Nanded, Maharashtra 431606, India}
\email{ushas073@gmail.com}
\begin{abstract} Hardy's theorem for the Riemann zeta-function $\zeta(s)$  says that it admits infinitely many complex zeros on the line  
$\Re({s}) = \frac{1}{2}$. In this note, we give a simple proof of this statement which, to the best of our knowledge, is new.
 \end{abstract}



 \maketitle


\section{Introduction}
\noindent
The Riemann Hypothesis (RH) is perhaps one of the most challenging problems in mathematics to this date. In 1859, Riemann  asserted that all the complex zeros of the Riemann zeta function $\zeta(s), s = \sigma + it $ lie on the \emph{critical line} $\sigma=1/2$. In other words, RH asserts that 
\begin{equation}
\zeta (s) \neq 0 \quad \textrm{in} \quad \sigma > 1/2.
\end{equation}
So far a proof or disproof of this statement has eluded the best minds. The sharpest known result in this direction is by I. M. Vinogradov \cite{V} who proved that there exists a positive constant $c$ such that 
\begin{equation}\label{vino}
\zeta (s) \neq 0 \quad \textrm{in} \quad \sigma > 1 - c(\log t)^{-2/3} (\log \log t)^{-1/3}.
\end{equation}
Riemann zeta-function was invented by Euler to study the distribution of primes. There is a deep connection between the zero-free region of $\zeta (s)$ and the error term in the prime number theorem (PNT) (see Chap. 12, \cite{ivic}).  Let $\pi (x)$ counts the number of primes upto $x$, then the prime number theorem says that there is an absolute constant $C>0$ such that
\[
\pi(x) = \int_{2}^{x} \frac{1}{\log u} du + O \left\{ x \exp \left(-C(\log x)^{3/5} (\log\log x)^{-1/5}  \right) \right\}.
\]
The above error term is a consequence of \eqref{vino}, which is the best known unconditional result whereas the error on RH  is $O(x^{{1/2}+\varepsilon})$.

It was G. H. Hardy who in 1914  \cite{hardy} showed that $\zeta(s)$ has infinitely many zeros on the critical line $ \sigma= 1/2 $.  This is generally known as a first step towards RH!  The method of Hardy paved way to establish analogous results for more general Dirichlet Series (see \cite{srini}). Titchmarsh's famous book on Riemann zeta-function (Chap. X, \cite{T}) discusses several methods of proving Hardy's theorem.

The basic idea behind the proof of Hardy's theorem is very simple. 
Suppose $f(x)\in C^{1}[a,b],$  and  $f(x)\neq 0 $ in $a\leqslant x \leqslant b $, then 
\begin{equation}\label{basic}
\vert\int_{a}^{b} f(x) dx \mid  =  \int_{a}^{b} \mid f(x) \mid dx.
\end{equation}
Therefore, if \eqref{basic} is violated for some $x$, then it shows that $f(x)$ has a zero of \emph{odd} order in $[a,b]$.

The proof of Hardy's theorem involves estimating a certain integral both from above and below and arrive at a contradiction by comparing these estimates. The estimation from below is straight forward for the Riemann zeta function and not too difficult in the case of a general Dirichlet series - thanks to a beautiful theorem of Ramachandra (see \cite{srini}). However, obtaining good upper bounds, which will lead to a contradiction, is a difficult problem!

In this article, we shall first discuss the classical approach of Hardy's theorem as given in  \S{10.5} of \cite{T}. Then in the next section we shall give an alternative proof of this theorem, which is simpler and seems to be new!

\medskip

\section{\textsc{Hardy's theorem for the Riemann zeta-function}}\label{H-T}
The Riemann zeta function, denoted by $\zeta(s),s=\sigma+it $, is defined by the series 
$$
\zeta(s)=\sum_{n=1}^{\infty} \frac{1}{n^{s}}  ,\quad \sigma > 1. 
$$            
It satisfies the following analytic properties:
\begin{enumerate}
\item[i)] Analytic Continuation: $(s-1)\zeta(s) $ is an entire function of order $1$.
\item [ii)] Functional Equation:
\begin{equation}\label{fe-zeta}
 \pi^{-s/2}\Gamma(s/2)\zeta(s) = \pi^{-(1-s)/2} \Gamma ((1-s)/2 )\zeta(1-s).
 \end{equation}
\item [iii)] Euler Product:
$$
\zeta(s)=\prod_{p}(1- {p^{-s}})^{-1}, \quad     \sigma >1,
$$
\end{enumerate}
\noindent
where the product is over all primes.

It is well known that the gamma function $\Gamma(s)$ is analytic in $\sigma >0$ and has no zeros. It has simple poles at $s = 0, -1, -2, -3, \cdots$.This fact, along with the functional equation \eqref{fe-zeta} implies that $\zeta(s)$ has zeros at $s=-2,-4,-6,...$.These zeros on the real line are called \emph{trivial} zeros of $\zeta(s)$.
From the theory of entire functions of finite order, it can be derived that all the complex zeros of $\zeta(s)$ have real part $0 \leq \sigma \leq 1$. It is also well known that $\zeta(1+it) \neq 0$ for $ t\neq 0$ and therefore $\zeta (it) \neq 0$ by \eqref{fe-zeta}. Further, it is an elementary exercise to show that $\zeta(s) \neq 0$ for $ 0 < s < 1$. Thus, all the complex zeros of $\zeta(s)$ lie in the vertical strip $0 < \sigma < 1$. These zeros are called \emph{non-trivial} zeros and the 
Riemann Hypothesis asserts that all the non-trivial zeros are, in fact, on the line $\sigma = 1/2$.

\medskip

Now, the theorem of Hardy for the Riemann zeta-function is the following:

\medskip

\noindent
\textbf{Hardy's Theorem.} \emph{$\zeta(s)$ has infinitely many zeros on the line $\sigma={1}/{2}.$}

\medskip

\noindent
\textbf{Proof.}
The functional equation \eqref{fe-zeta} can be rewritten as 
\begin{equation}\label{fe-zeta-2}
\zeta(s) = \chi(s) \zeta(1-s), 
\end{equation}
where $\chi(s) = {\pi^{s-1/2} \Gamma ( \tfrac{1}{2}(1-s))}/{\Gamma ( \tfrac{1}{2} s )}.$

Observe that $\chi(s)\chi(1-s) = 1$ or $ \chi(\tfrac{1}{2}+it) \chi(\tfrac{1}{2}-it) = 1, $
i.e., $\vert\chi(\frac{1}{2}+it)\vert = 1.$ From \eqref{fe-zeta-2}, we get 
$ \chi(\tfrac{1}{2}+it)^{-1/2}  \zeta(\tfrac{1}{2}+it) = \chi(\tfrac{1}{2}-it)^{1/2}  \zeta(\tfrac{1}{2}-it).$

Thus the function $$Z(t) :=  \chi(\tfrac{1}{2}+it)^{-1/2}  \zeta(\tfrac{1}{2}+it)$$ is \emph{real} for real values of $t\neq 0$. Moreover, $ \vert Z(t) \vert  = \vert \zeta(\tfrac{1}{2} + it) \vert $. The function $Z$ is popularly called Hardy's $Z$-function.  From the above discussion, it is clear that the zeros of $Z(t)$ correspond to the zeros of $\zeta (s) $ on the critical line. We shall apply the basic idea elucidated in the beginning by taking $ f = Z $ in \eqref{basic}.

To begin with,  we shall assume that  $Z(t)$ has finitely many zeros or no zeros, then we can find a real  number $T_{0}>0 $ sufficiently large such that for all $ T \geq T_{0}$, we have 
\begin{equation}\label{basic-2}
\vert I \vert := \vert \int_{T}^{2T} Z(t)dt \vert  = \int_{T}^{2T}\vert Z(t) \vert dt .
\end{equation}

Now, we need to estimate the integrals in \eqref{basic-2} both from below and above and arrive at a contradiction to our assumption. 

\noindent 
\textbf{Estimation from below:} This is the easy step. We have
\begin{equation}\label{lower-b-1}
\int_{T}^{2T}\vert Z(t) \vert dt = \int_{T}^{2T}\vert \zeta(\tfrac{1}{2}+it) \vert dt \geq  \vert \int_{T}^{2T}\zeta(\tfrac{1}{2}+it) dt \vert. 
\end{equation}
Now, by a change of variable we get
$$ 
i \int_{T}^{2T}\zeta(\tfrac{1}{2}+it) dt = \int_{{1}/{2}+iT}^{{1}/{2}+2iT}\zeta(s)ds.
$$
Moving the line of integration to $\sigma = 2$ and using Cauchy's theorem, we obtain
\begin{equation}\label{lower--2}
\int_{{1}/{2}+iT}^{{1}/{2}+2iT}\zeta(s)ds
 =(\int_{{1}/{2}+iT}^{2+iT}+ \int_{2+iT}^{2+2iT}+\int_{2+2iT}^{{1}/{2}+2iT})\zeta(s)ds. 
\end{equation}

\noindent
Now, from the convexity bound \eqref{zeta-b}, it follows that
 \begin{equation}\label{lower-b-3}
 \int_{{1}/{2}+iT}^{2+iT}\zeta(s) ds = \int_{{1}/{2}}^{2}\zeta(\sigma + iT) d\sigma  = O \left( T^{1/4 +\varepsilon}   \right).
 \end{equation}

\noindent
On the other hand,
\begin{eqnarray}\label{lower-b-4}
\int_{2+iT}^{2+2iT}\zeta(s)ds  & = & \int_{2+iT}^{2+2iT}\left(1+ \sum_{n=2}^{\infty} {n^{-s} } \right) ds \nonumber \\
& =  & iT +  i \sum_{n=2}^{\infty} {n^{-2}}\int_{T}^{2T} e^{-it \log n} dt \nonumber\\
& = &  iT + O\left( \sum_{n=2}^{\infty} \tfrac{1}{{n^2}\log n} \right) = iT + O(1).
\end{eqnarray}
\noindent
Therefore, from \eqref{lower-b-1}, \eqref{lower--2}, \eqref{lower-b-3} and  \eqref{lower-b-4}, we obtain
\begin{equation}\label{lb} 
\int_{T}^{2T}\vert Z(t) \vert dt  > A T. 
\end{equation}
where $A > 0$ is an absolute constant.
\medskip

\noindent
The main problem, as remarked before, is to find an upper bound of the type
\begin{equation}\label{ub-0}
 \vert I \vert \ll T^{\alpha} .  
\end{equation}
If  $0 < \alpha < 1$, this will contradict \eqref{lb} and thus establishes the theorem.

\medskip

\noindent
\textbf{Remark.} As $ \vert \chi(\tfrac{1}{2}+it) \vert =1$, taking trivial estimate we see that
\begin{eqnarray*}
\vert  \int_{T}^{2T}Z(t)dt  \vert & = & \vert \int _{T}^{2T} \chi(\tfrac{1}{2}+it)^{-1/2}\zeta(\tfrac{1}{2}+it) dt \vert \\
{} & \leq &  \int _{T}^{2T} \vert \zeta(\tfrac{1}{2}+it)  \vert dt. 
\end{eqnarray*}
Even on Lindel\"{o}f hypothesis (which is a consequence of RH and which says that $ \zeta(\sigma +it) \ll t^{\varepsilon}$ for $\sigma \geq \tfrac{1}{2}$) the above is only $ O ( T^{1+\varepsilon} )$, not \emph{enough} to contradict \eqref{lb}.

\medskip

\noindent
Thus we need to estimate $ \vert I \vert $ nontrivially, which we do below.

 \medskip
 
\noindent
\textbf{Estimation from above :}

\medskip

Our aim is to estimate the integral
\begin{equation*}
\vert I \vert  = \vert \int_{1/2 + iT}^{1/2 + 2iT} \chi(s)^{-1/2} \zeta(s) ds \vert.
\end{equation*}
We move the line of integration to the right of the line $\sigma = 1$ so that $\zeta(s)$ can be replaced by its series representation and then evaluate the integral using Cauchy's theorem. Accordingly, let us consider the rectangle $\mathcal{R}$ with the vertices at ${1}/{2}+iT, 1+\delta+iT, 1+\delta+2iT, {1}/{2}+2iT$ for some $\delta > 0$ to be chosen later.

\medskip

As $\chi(s) \neq 0$ in $\sigma > 0, \  \chi(s)^{-1/2}$ is single-valued analytic function in $\sigma \geq 1/2$, therefore, by Cauchy's Theorem
\[
\frac{1}{2 \pi i}\int_{\mathcal{R}}\chi(s)^{-1/2} \zeta(s)ds = 0. 
\]

\medskip

\noindent
Thus,
\begin{align}\label{cauchy-1}
\frac{1}{2\pi i}\int_{{1}/{2}+iT}^{{1}/{2}+2iT} \chi(s)^{-1/2} \zeta(s)ds 
 & = \frac{1}{2\pi i}\left(\int_{{1}/{2}+iT}^{1+\delta +iT} + \int_{1+\delta+iT}^{1+\delta +2iT} + \int_{1+\delta+2iT}^{{1}/{2}+2iT}\right)\chi(s)^{-1/2}\zeta(s) ds     \nonumber \\
 & = I_{1}+I_{2}+I_{3}       \text{(say)}. 
\end{align}
As $I_{1}$ and $I_{3}$ are of the same order of magnitude, it is  enough to estimate $I_{1}$ only.

We require growth estimates for the integrand in the region $ \sigma \geq 1/2$.  The well-known Stirling's formula for $\Gamma(s)$ says that  in any fixed strip $\alpha \leq \sigma \leq \beta$, we have, 
$$\Gamma(\sigma+it) = t^{\sigma+it-\frac{1}{2}} e^{\frac{1}{2} \pi t -it + \frac{1}{2} i \pi (\sigma -\frac{1}{2})}(2\pi)^{\frac{1}{2}} \left( 1 + O (t^{-1}) \right),$$ as $t \to \infty$. This gives
\begin{equation}\label{chi-1}
\chi(s) = \left( {t}/{2\pi}\right )^{{1}/{2}-(\sigma+it)} e^{i(t+{\pi}/{4})} \left(1+O(t^{-1})\right),
\end{equation}
as  $ t \rightarrow \infty$.

\medskip

\noindent
Therefore,
\begin{equation}\label{chi-2}
\chi(s)^{-1/2} = O\left(t^{\sigma/2 - 1/4}\right),        \quad         0\leq \sigma\leq 1+\delta. 
\end{equation}
We also require an upper bound for $\zeta (s)$ in the critical strip which is uniform in $\sigma$. This is obtained by using convexity principle. First, one shows that $\zeta(s) = O(\log t)$ for $\sigma \geq 1$. Then \eqref{fe-zeta-2}  combined with \eqref{chi-1} gives $ \zeta(0+it) = O(t^{1/2} \log t)$. Finally the  Phragmen-Lindel\"of convexity principle implies 
\begin{equation}\label{zeta-b}
\zeta(s) = O \left( t^{(1-\sigma)/2} (\log t)^5 \right), \mbox{uniformly in} \     0 \leq \sigma \leq 1.
\end{equation}

\medskip

Now, combining \eqref{chi-2} and \eqref{zeta-b}, we obtain
\[
\chi(s)^{-1/2}\zeta(s) = \begin{cases}
O\left( t^{1/4} (\log t)^5 \right), & 0 \leq \sigma \leq 1 {} \\
O( t^{1/4 + {\delta}/2}),  & 1 \leq \sigma \leq 1+\delta {}. 
\end{cases}
\end{equation*}

Thus,
\begin{equation}\label{ub-h}
\vert I_1\vert = \vert \int_{1/2}^{1+\delta} \chi(\sigma +iT)^{-1/2}\zeta(\sigma+iT)  \ d\sigma \vert =  O\left(T^{{1}/{4}+\delta/2}\right).
\end{equation}

Now comes the difficult part, the estimation of $I_2$, which we carry out below.

\medskip

\noindent
We have 
\begin{align*}
 I_{2}  = &  \frac{1}{2\pi i}\int_{1+\delta+iT}^{1+\delta+2iT}\chi(s)^{-1/2} \zeta(s) ds \\
  = & \frac{1}{2\pi} \int_{T}^{2T}\chi(1+\delta+it)^{-1/2} \zeta(1+\delta+it)) dt. 
\end{align*}
By Stirling's formula \eqref{chi-1}, we have
\[
 \chi(1+\delta+it)^{-1/2} = \left({\tfrac{t}{2\pi}}\right)^{{1}/{4}+{\delta}/{2}+{it}/{2}} e^{-i ({t}/{2}+{\pi}/{8})} \left\{1+O({1}/{t})\right\}.
\]
Therefore, the above integral is
\begin{align*}
 & = \sum_{n=1}^{\infty}\tfrac{1}{n^{1+\delta}}\int_{T}^{2T}\left({\tfrac{t}{2\pi}}\right)^{{1}/{4}+{\delta}/{2}+{it}/{2}} e^{-i ({t}/{2}+{\pi}/{8})} n^{-it} \left\{1+O({1}/{t})\right\} dt.
 \end{align*}
The $O$-term is $ O(\int_{T}^{2T}t^{{-3/4}+\delta/2} dt) = O(T^{1/4 + \delta/2}) = o(T) \quad \mbox{provided} \quad \delta < 3/2. $

The main term is a constant multiple of 
 \begin{align}\label{ub-sum}
 {} & \sum_{n=1}^{\infty} \tfrac{1}{n^{1+\delta}} \int_{T}^{2T} \left({\tfrac{t}{2\pi}}\right)^{{1}/{4}+{\delta}/{2}+{it}/{2}} e^{-{it}/{2}-it \log n} dt  \nonumber \\
 {} & = \sum_{n=1}^{\infty} \tfrac{1}{n^{1+\delta}} \int_{T}^{2T} \left({\tfrac{t}{2\pi}}\right)^{{1}/{4}+{\delta}/{2}} e^{\frac{it}{2} \log ({t}/{2\pi en^{2})}}  dt.
 \end{align}

 Observe that a trivial estimation of the integral in \eqref{ub-sum} will not yield the desired upper bound, thus we need to exploit the possible cancellations arising from the exponential term in the integrand and hope that the resulting estimate will be \emph{small}!
 
The following results are  standard tools in theory of exponential integral techniques (see page 71--72, \cite{T}) . 
\begin{lemma}\label{l-0}
Let $F(x)$ be a real-valued differentiable function such that $F'(x)$ is monotonic, and $F'(x) \geq m > 0$, or $F'(x) \leq -m < 0$, throughout the interval $[a, b]$. Then
\[
\vert \int_{a}^{b} e^{i F(x)} dx \vert \leq \frac{4}{m}.
\]
\end{lemma}

\begin{lemma}\label{l-1}
Let $F(x)$be twice differentiable function and $F''(x) \geq r  > 0$ or $F''(x) \leq -r < 0$,  throughout the interval $[a,b]$ and let ${G(x)}/{F'(x)}$ be monotonic with $\vert G(x)\vert \leq M$, then 
\[
\vert \int_{a}^{b} G(x) e^{iF(x)}dx \vert \leq \frac{8M}{\sqrt{r}}.
\]
\end{lemma}

Now, taking $F_n(t) =  \frac{t}{2} \log( \tfrac{t}{2\pi en^{2}}) $, we see that  $  {F''}_{n}(t) =1/2t \geq 1/4T, T\leq t \leq 2T$.

Thus from the Lemma \ref{l-1} we obtain
\begin{equation}\label{ub-2}
\int_{T}^{2T} {\left(\tfrac{t}{2\pi}\right)}^{{1}/{4}+{\delta}/{2}} e^{ \frac{it}{2} \log ({t}/{2\pi en^{2})}} dt = O \left( T^{{1}/{4}+{\delta}/{2}}T^{{1}/{2}} \right) = O \left( T^{{3}/{4}+{\delta}/{2}} \right). 
\end{equation}

Now combining \eqref{cauchy-1}, \eqref{ub-h}, \eqref{ub-sum} and \eqref{ub-2}, we obtain
\begin{equation}\label{ub-final}
 \vert\int_{T}^{2T} Z(t) dt \vert = O(T^{{3}/{4}+\delta/2}) = o(T) \ \text{provided} \ \delta < 1/2.
 \end{equation}

Therefore, the desired contradiction is arrived from  comparing \eqref{lb} and \eqref{ub-final}, by choosing $\delta <1/2.$

\medskip

\noindent
\textbf{Remark.} Another method of estimating the upper bound for $\vert I \vert$ uses the following weak version of the well-known Riemann-Siegel formula
\begin{equation}\label{R-S}
Z(t) = 2 \sum_{ n\leq \sqrt{t/(2\pi)}} n^{-1/2} \cos \left( t \log \frac{ \sqrt{t/(2\pi)}}{n} - \frac{t}{2} - \frac{\pi}{8}  \right) + O\left( t^{-1/4} \right).
\end{equation}
The integral $I$ is then evaluated using Lemma {\ref{l-1}} to obtain $I = O(T^{3/4})$ (see Chap. 2 of \cite{ivic-z} for details).  

\section{A new approach using first approximation formula for $\zeta(s)$}

In this method, the integral
\begin{equation}\label{ub-3}
 I = \int_{T}^{2T}\chi(\tfrac{1}{2}+it)^{-1/2}\zeta(\tfrac{1}{2}+it) dt
 \end{equation} 
 is evaluated by replacing the zeta function with a Dirichlet polynomial having a \emph{good} error term. As mentioned in the remark above, a weak version of Riemann-Siegel formula can be used for the upper bound estimation.  
 We will use the  well-known first order approximation of $\zeta (s)$ (see page 77 of \cite{T}) to estimate $I$ both from below and above. The estimation from above seems to be new.
\begin{lemma} We have
\begin{equation}\label{zeta-app-1}
 \zeta(s) = \sum_{n \leq x}{n^{-s}} + \frac{x^{1-s}}{s-1} + O(x^{-\sigma}),  
\end{equation}
 uniformly for $\sigma \geq \sigma_0 >0,  \vert t \vert  < {2 \pi x}/{C} $
 when $C$ is a given constant $> 1$.
 \end{lemma}

Taking $x = C T/\pi$ in the above Lemma and noting that $ T \leq t \leq 2T$, we get
 \begin{equation}\label{first-e-1}
 \zeta(\tfrac{1}{2} +it) = \sum_{n \leq CT/\pi}{n^{-\tfrac{1}{2} - it}} + O \left(\frac{T^{1/2}}{\vert t \vert} \right)+ O(T^{-1/2}).  
\end{equation}

For the lower bound, recall that we have
\begin{equation}\label{lower-b-1-1}
\int_{T}^{2T}\vert Z(t) \vert dt = \int_{T}^{2T}\vert \zeta(\tfrac{1}{2}+it) \vert dt \geq  \vert \int_{T}^{2T}\zeta(\tfrac{1}{2}+it) dt \vert. 
\end{equation}
Therefore, from \eqref{first-e-1} we have
\begin{eqnarray}
\int_{T}^{2T} \zeta(\tfrac{1}{2}+it) dt &=& \int_{T}^{2T} \{ 1+\sum_{2\leq n \leq {CT}/{\pi}} {n^{-\tfrac{1}{2} - it}} \} dt + O(T^{1/2})\\
& = & T + O (T^{1/2})  > A T,
\end{eqnarray}
where $A > 0$ is an absolute constant.

Therefore, we have
\begin{equation}\label{lower-bound}
\int_{T}^{2T}\vert Z(t) \vert dt \geq AT.
\end{equation}

\medskip
\noindent
On the other hand, for estimation of upper bound, we have from \eqref{first-e-1},
\begin{equation}\label{first-e-2}
Z(t) =  \chi(\tfrac{1}{2}+it)^{-1/2} \sum_{n\leq CT/\pi} {n^{-\tfrac{1}{2} - it}} + O \left(\frac{T^{1/2}}{\vert t \vert} \right)+ O(T^{-1/2}).  
\end{equation}

\noindent
The $O$-terms contributes $O ( T^{1/2} ) $ to the integral $I$.
Thus, by Stirling's formula for $\chi(s)$, we have
\begin{equation}
I  =  \int_{T}^{2T}  \sum_{n\leq CT/\pi} {n^{-\tfrac{1}{2}-it}} \left(\tfrac{t}{2\pi}\right)^{{it}/{2}} e^{-i({t}/{2}+{\pi}/{8}) } \left\{  1+O\left(\tfrac{1}{t}\right) \right\} dt + O(T^{1/2}).  
\end{equation}
Note that the error term in the above expression is $O ( T^{1/2} )$.  Thus the main term is a constant multiple of
\begin{equation}
\sum_{n\leq CT/\pi}\tfrac{1}{n^{{1}/{2}}}\int_{T}^{2T} e^{i F(n,t) }  dt, 
\end{equation}
where 
$ F(n,t) = \tfrac{t}{2} \log \left(\tfrac{t}{2\pi en^{2}}\right)$. Observe that 
\begin{equation}\label{d-bound}
F'(n,t) =  \tfrac{1}{2} \log \left( \tfrac{t}{2\pi n^{2}} \right), \  F''(n,t)= \tfrac{1}{2t} > \tfrac{1}{4T}.
\end{equation}
As $ F'(n,t) = 0$ when $ t = 2\pi n^2 \, \text{or when} \,  n = \sqrt{t / 2\pi}$,  we break the sum over $n$ into two parts viz., $ 1 \leq n \leq 3\sqrt{T/\pi} $ and $ 3\sqrt{T/\pi} < n \leq CT/\pi$. The integral in the first sum is evaluated using Lemma \ref{l-1} and the integral in the second sum is evaluated using Lemma \ref{l-0}.

Thus from Lemma \ref{l-1} and \eqref{d-bound} we obtain
\begin{equation}\label{sum-1}
\vert \sum_{n\leq 3\sqrt{T/\pi}}\frac{1}{n^{{1}/{2}}}\int_{T}^{2T}  e^{i F(n,t)} dt\vert = O\left( \sum_{n\leq{3\sqrt{T/\pi}}} \frac{1}{n^{1/2}} T^{1/2}\right) = O\left(T^{3/4} \right).
\end{equation}

To evaluate the second sum, we use the elementary inequality $ \log \frac{n}{m} \geq \frac{n-m}{n}$ for $ n> m$. 

This gives
\[
 \log \left( \tfrac{2\pi n^{2}}{t}\right) \geq \left( \tfrac{2\pi n^{2}-t}{2 \pi n^2}\right) \geq \tfrac{8}{9},
\]
as $ T \leq t \leq 2T$ and $ 3\sqrt{T/\pi} < n \leq CT/\pi$. Thus
$  F'(n,t) = \tfrac{1}{2} \log \left( \tfrac{t}{2\pi n^{2}} \right) \leq -4/9$. Now applying Lemma \ref{l-0}, we get
\begin{equation}
\vert \sum_{3\sqrt{T/\pi} < n \leq CT/\pi}\frac{1}{n^{{1}/{2}}}\int_{T}^{2T} e^{i F(n,t)} dt\vert = O\left(\sum_{3\sqrt{T/\pi} < n \leq CT/\pi} \frac{1}{n^{1/2}} \right) = O\left(T^{1/2} \right).
\end{equation}
Thus, it follows that 
\begin{equation}\label{upper-bound}
\vert\int_{T}^{2T} Z(t) dt \vert = O(T^{3/4}).
\end{equation}

\medskip

\noindent
Equation \eqref{lower-bound} and \eqref{upper-bound} leads to a contradiction to our assumption and thus establishes Hardy's theorem.

\medskip

\noindent
\textbf{Acknowledgements.} I wish to thank Prof. K. Srinivas  for many useful discussions. I also thank Prof. A. Ivi\'c for his advise and  encouragement. I express my gratitude to the Institute of Mathematical Sciences for providing congenial environment which enabled me to carry out this project. I sincerely thank Prof. C. S. Aravinda for the prompt handling of the manuscript and for suggesting  some changes in the layout of article.

\end{document}